\documentclass[11pt]{amsart}%
\usepackage{amsmath}%
\usepackage{amssymb}%
\usepackage{amscd}%
\usepackage{color}%
\usepackage[all]{xy}

\renewcommand{\Bbb}{\mathbb}  

\newcommand{\Na}{{\Bbb{N}}} 
\newcommand{\Q}{{\Bbb{Q}}}  
\newcommand{\R}{{\Bbb{R}}}  
\newcommand{\Z}{{\Bbb{Z}}}  
\newcommand{\G}{{\Bbb{G}}}  
\renewcommand{\H}{{\Bbb{H}}}

\newcommand{\F}{{\Bbb{F}}}  

\newcommand{\Oh}{\mathcal{O}}

\newcommand{\tensor}{\otimes}   
\newcommand{\Mor}{\operatorname{Mor}}
\newcommand{\Aut}{\operatorname{Aut}}
\newcommand{\prolim}{\varprojlim}
\newcommand{\Der}{\mathrm{Der}}
\newcommand{\ord}{\mathrm{ord}}

\newcommand{\g}{\mathfrak{g}}

\newcommand{\Log}{\operatorname{Log}}
\newcommand{\GG}{\mathcal{G}}
\newcommand{\II}{\mathcal{I}}
\newcommand{\LL}{\mathcal{L}}
\newcommand{\gr}{\mathrm{gr}}
\newcommand{\Sat}{\mathrm{Sat}}
\newcommand{\la}{\mathrm{la}}

\def\wh#1{\widehat{#1}}
\theoremstyle{plain}
    \newtheorem{theorem}{Theorem}[subsection]
    \newtheorem{thm}[theorem]{Theorem}
    \newtheorem{appthm}{Theorem}[section] 
    
    \newtheorem{proposition}[theorem]{Proposition}
    \newtheorem{prop}[theorem]{Proposition}
    \newtheorem{appprop}[appthm]{Proposition} 
    \newtheorem{lemma}[theorem]{Lemma}
    
     \newtheorem{cor}[theorem]{Corollary}
\theoremstyle{definition}
    \newtheorem{definition}[theorem]{Definition}
    \newtheorem{defn}[theorem]{Definition}
    
    \newtheorem{example}[theorem]{Example}
    
    \newtheorem{lemmadef}[theorem]{Definition and Lemma}
    \newtheorem{remark}[theorem]{Remark}
    \newtheorem{rem}[theorem]{Remark}
 
\def\Alphabet{A,B,C,D,E,F,G,H,I,J,K,L,M,N,O,P,Q,R,S,T,U,V,W,X,Y,Z}
\def\alphabet{a,b,c,d,e,f,g,h,i,j,k,l,m,n,o,p,q,r,s,t,u,v,w,x,y,z}
\def\endpiece{xxx}
\def\makeAlphabet[#1]{\expandafter\makeA#1,xxx,}
\def\makealphabet[#1]{\expandafter\makea#1,xxx,}
\def\makeA#1,{\def\temp{#1}\ifx\temp\endpiece\else%
\mkbb{#1}\mkfrak{#1}\mkbf{#1}\mkcal{#1}\expandafter\makeA\fi}%
\def\makea#1,{\def\temp{#1}\ifx\temp\endpiece\else\mkfrak{#1}\mkbf{#1}\expandafter\makea\fi}%
\def\mkbb#1{\expandafter\def\csname bb#1\endcsname{\mathbb{#1}}}
\def\mkfrak#1{\expandafter\def\csname fr#1\endcsname{\mathfrak{#1}}}
\def\mkbf#1{\expandafter\def\csname b#1\endcsname{\mathbf{#1}}}
\def\mkcal#1{\expandafter\def\csname c#1\endcsname{\mathcal{#1}}}
\def\makeop[#1]{\xmakeop#1,xxx,}
\def\mkop#1{\expandafter\def\csname #1\endcsname{{\mathrm{#1}}}} %
\def\xmakeop#1,{\def\temp{#1}\ifx\temp\endpiece\else\mkop{#1}\expandafter\xmakeop\fi}%
\makeAlphabet[\Alphabet]
\makealphabet[\alphabet]
\makeop[Alt,End,Ext,Hom,Sym,Tor]
\makeop[CH,Pic,div]
\makeop[Ker,Im,Coker,Coim,id,pr,isom]
\makeop[Spec,Spf,Spm]
\makeop[Re,Im]%
\makeop[dR,Nis,crys,rig,syn,tor,Cont]
\makeop[Gal,Res,ab]
\makeop[pol]
\makeop[Var,an,Isoc,univ]%
\makeop[Gl,Sl]
\makeop[Lie]
\makeop[tay,tr]
\makeop[Sat]

\def\isom{\cong}

\def\prolim{\varprojlim}

\begin{document}
\title[Lazard isomorphism]{Some complements to the Lazard isomorphism}   
\author{Annette Huber, Guido Kings and Niko Naumann}
\date{24 February, 2010 }  
\begin{abstract} Lazard showed in his
seminal work \cite{L} that for rational coefficients 
continuous group cohomology of $p$-adic Lie-groups 
is isomorphic to Lie-algebra cohomology.
We refine this result in two directions:
firstly we extend his isomorphism under certain conditions to integral coefficients and
secondly, we show that for algebraic groups over finite extensions $K/\Q_p$ his isomorphism can
be generalized to $K$-analytic cochains and $K$-Lie algebra cohomology.
\end{abstract}

\maketitle
\tableofcontents

\section{Introduction}

One of the main results of Lazard's magnum opus \cite{L} on $p$-adic Lie-groups
is a comparison isomorphism between continuous group cohomology,
analytic group
cohomology and Lie-algebra cohomology. This comparison isomorphism 
is an important tool in the cohomological study of Galois representations
in arithmetic geometry. More recently, it also appeared in topology and
homotopy theory in connection with the formal groups associated to
cohomology theories and in particular with topological modular forms.

Lazard's comparison theorem holds for $\bbQ_p$-vector spaces and the
isomorphism between continuous cohomology and Lie-algebra cohomology
is obtained from a difficult isomorphism between the saturated group ring
and the saturated universal enveloping algebra. For some applications
(e.g. the connection with the Bloch-Kato exponential map  in \cite{HK})
it is important to have a version for  integral coefficients and 
a better understanding of the map between the cohomology theories.

In this paper we extend and complement the comparison isomorphism in
two directions:

The first result is an integral version of the isomorphism
assuming technical conditions (Theorem
\ref{intlazard}). For uniform pro-$p$-groups one gets a clean result
with only a mild condition on the module (Theorem \ref{uniform}).
To our knowledge and with the notable exception of Totaro's work
\cite{T}, this is the first progress on a problem which Lazard wrote
``reste \`a faire'' more than forty years ago \cite[Introduction 7,C)]{L}.

The second result concerns the definition of the isomorphism in the 
case of smooth group schemes. Here one can define directly a map from analytic
group cohomology to Lie-algebra cohomology with constant coefficients 
by differentiating cochains (see Definition \ref{analLaz}). We showed in \cite{HK}
(see  Proposition \ref{lazard}) that the resulting map is Lazard's comparison isomorphism modulo the identification of continuous cohomology with analytic
cohomology. Serre mentioned to us that
this was clear to him at the time Lazard's paper was
written, however, it was not included in the published results. Unfortunately, we were so
far not able to use this simple map to obtain an independent proof of Lazard's comparison result.

The advantage of this description of the map is not only its simplicity but also that it carries over to
$K$-Lie-groups for finite extensions $K/\Q_p$. In Theorem \ref{isom} we prove that this map is also
an isomorphism in the case of
$K$-Lie-groups attached to smooth group schemes with connected generic
fiber over the integers of $K$. This generalizes results in \cite{HK}
for $\mathrm{GL}_n$ and complements the result of Lazard, who only treats $\bbQ_p$-analytic groups.

The paper is organized as follows. In Section \ref{review}
we give a quick tour through the notions from \cite{L} that we need.
We hope that this section also proves to be  a useful overview 
of the central notions and results in \cite{L}.
In Section \ref{intsection} we prove our integral
refinement of Lazard's isomorphism. Finally Section
\ref{groupsection} considers the isomorphism over a general base in
the case of group schemes.

{\em Acknowledgments. We would like to thank B.~Totaro for numerous insightful
remarks on a preliminary draft of this paper. The third author thanks
H.-W.~Henn for discussions about Example \ref{morava}. Finally, the authors
thank T. Weigel for an insightful remark explained in the Appendix.}

\section{Review of some results by Lazard}\label{review}
In this section we recall the basic notions about groups and
group rings we need to formulate our main results. As we proceed we illustrate the
main notions with the example of separated smooth  group schemes. We hope
that this section helps to guide through the long and difficult paper by Lazard.

\subsection{Saturated groups}

\begin{definition}\cite[II 1.1, II 1.2.10, III 2.1.2]{L}\label{filtrationdefn}
A {\em filtration} on a group $G$ is a map
\[ \omega: G\to \R^*_+\cup\{\infty\}\]
such that 
\begin{itemize}
	\item[1)] For all $x,y\in G$, $\omega(xy^{-1})\ge \inf\{\omega(x),\omega(y)\}$
	\item[2)] For all $x,y\in G$, $\omega(x^{-1}y^{-1}xy)\ge \omega(x)+\omega(y)$.
\end{itemize}
$G$ is called {\em $p$-filtered} if in addition for all $x\in G$
$$
	\omega(x^p)\geq \inf\{\omega(x)+1,p\omega(x)\}.
$$
$G$ is called {\em $p$-valued} if $\omega$ satisfies 
\begin{itemize}
	\item[3)] $\omega(x)<\infty$ for $x\neq e$
	\item[4)]  $\omega(x)>(p-1)^{-1}$ for all $x\in G$
	\item[5)]   $\omega(x^p)=\omega(x)+1$ for all $x\in G$.
\end{itemize}
$G$ is called {\em $p$-divisible} if it is $p$-valued and 
\begin{itemize}
	\item[6)] For all $x\in G$ with $\omega(x)>1+\frac{1}{p-1}$ there exists $y\in G$ with
	$y^p=x$.
\end{itemize}
Finally, a $p$-divisible group $G$ is \emph{saturated}, if
\begin{itemize}
	\item[7)] $G$ is complete for the topology defined by the filtration. 
\end{itemize}
\end{definition}
Note that a filtration satisfying the conditions of a $p$-valuation is
automatically $p$-filtered.

Recall that a pro-$p$-group is the inverse limit of finite $p$-groups.
This is the case we are going to work with.
 
\begin{proposition}\cite[II 2.1.3]{L}\label{pro-p-crit} A $p$-filtered group is a pro-$p$-group if and only 
if it is compact.
\end{proposition}

We denote by $\F_p[\epsilon]$ the polynomial ring with generator $\epsilon$ in
degree $1$.

\begin{lemmadef}\cite[II 1.1., III 2.1.1, III 2.1.3]{L}
Let $(G,\omega)$ be filtered. 
\begin{enumerate}
\item
For every $\nu\in\R^*_+$
\begin{equation*} G_\nu:=\{ x\in\G\, |\, \omega(x)\ge\nu\}, G_\nu^+:=\{ x\in G\, |\, \omega(x)>\nu\}\subseteq G\end{equation*}
are normal subgroups.
\item
\begin{equation*} \gr(G):=\bigoplus_{\nu\in \R^*_+} G_\nu/G_\nu^+,\end{equation*}
$\gr(G)$ is a graded Lie-algebra over $\F_p$, the Lie-bracket being induced by the
commutator in $G$.
\item If $(G,\omega)$ is $p$-valued, then $\gr(G)$
is even a graded $\F_p[\epsilon]$-Lie-algebra, the action of $\epsilon$ being induced by $x\mapsto x^p$ ($x\in G_\nu$, $x^p\in G_{\nu+1}$)
\item In this case, $\gr(G)$  is free as a graded
$\F_p[\epsilon]$-module.
The {\em rank of $G$} is by definition the rank of the $\F_p[\epsilon]$-module
$\gr(G)$.
\end{enumerate}
\end{lemmadef}

\begin{example} \cite[V 2.2.1]{L}
Let $(G,\omega)$ be a complete $p$-valued group of finite rank $d$.

If $\{x_i\}_{i=1,\ldots,d}\subseteq G$ are representatives of an ordered
basis of the $\F_p[\epsilon]$-module $\gr(G)$, then every $y\in G$ is uniquely
an ordered product $y=\prod\limits_{i=1}^d x_i^{\lambda_i}$  with 
$\lambda_i\in\Z_p$
and 
\[ \omega(y)=\inf\limits_i (\omega(x_i)+v(\lambda_i)),\]
where the valuation on $\Z_p$ is normalized by $v(p)=1$.
$(G,\omega)$ has rank $d$.
\end{example}
\begin{defn}\cite[V 2.2.1, V 2.2.7]{L}
\begin{enumerate}
\item In the situation of the example,
the family $\{x_i\}_{i=1,\ldots,d}$ is called an {\em ordered basis of $G$}.
\item The $p$-valued group $(G,\omega)$ is called {\em equi-$p$-valued} if there exists an
ordered basis $\{ x_i\}$ as above such that 
\[ \omega(x_i)=\omega(x_j)\mbox{ for all }1\leq i,j\leq d.\]
\end{enumerate}
\end{defn}

\subsection{Serre's standard groups as examples}\label{standardsection}
Let $E$ be a finite extension of $\Q_p$ with ring of integers $R$.
Let $\frm$ be the maximal ideal of $R$. Then
$E$ is a discretely valued field. 
We normalize its
valuation $v$ by $v(p)=1$. Let $e$ be the ramification index of $E/\Q_p$.

Any formal group law $F(X,Y)$ in $n$ variables over $R$ defines
a group structure $G$ on $\frm^n$. These are the \emph{standard groups} as defined by Serre.

Define for $(x_1,\ldots,x_n)\in \frm^n$
$$
	\omega(x_1,\ldots,x_n):=\inf_i\{v(x_i)\}.
$$
Then we have for $\lambda\ge 0$
$$
	G_\lambda:=\{x\in G:\omega(x)\ge \lambda\}.
$$
\begin{proposition} \cite[LG 4.23 Theorem 1 and Corollary and 4.25 Corollary 2 of Theorem 2]{Serre1}
	$G$ is a pro-$p$-group.
	For any $\lambda\ge 0$ the group $G_\lambda$ is a normal subgroup of $G$. Moreover,
	the map $\omega$ defines a filtration on $G$ in the sense of definition \ref{filtrationdefn},1)-3).
\end{proposition} 
One can show that a suitable open subgroup of the standard group
 $G$ is saturated.\\
Let $\rho$ denote the smallest integer which is larger than $\frac{e}{p-1}$.

\begin{lemma}\label{saturationlemma} 
Let $E$ and $G$ be as above. Then the subgroup $(H,\omega)$,
$$
	H:=\left\{x\in G:\omega(x)>\frac{1}{p-1}\right\}=G_\frac{\rho}{e}
$$
is saturated and of finite rank.
It is equi-$p$-valued if and only if  $e=1$. 
\end{lemma}
\begin{proof}
Note first that according to \cite[LG 4.21, Corollary]{Serre1} the power series
$f_p$ which defines the $p$-power map is of the form \[ f_p(X)=p(X+\varphi(X))+\psi(X)\]
with $\ord(\varphi)\ge 2$ and $\ord(\psi)\ge p$. It follows that
$$
	\omega(x^p)\ge \inf\{\omega(x)+1,p\omega(x)\},\, x\in G
$$
because if $x$ has coordinates $x_i$, then $x^p$ has coordinates
\[ f_{p,i}(x_1,\ldots,x_n)=px_i+p\varphi_i(x)+\psi_i(x)\]
and the valuations of the summands are bounded below by $1+\omega(x)$,
$1+2\omega(x)>1+\omega(x)$ and $p\omega(x)$, respectively.\\
As $\omega(x)>\frac{1}{p-1}$ is equivalent to $\omega(x)+1<p\omega(x)$,
this implies on $H$ 
\[ \omega(x^p)\ge \omega(x)+1.\]
On the other hand, let  $x_i$ be a coordinate of $x$ with $\omega(x)=v(x_i)$. Then
\[ \omega(x^p)\leq v(px_i+p\phi_i(x)+\psi_i(x))=1+\omega(x_i)=1+\omega(x)\]
and hence
\[ \omega(x^p)=\omega(x)+1\]
for all $x\in H$. This shows that  $(H,\omega)$ is $p$-valued. 
To see that $H$ is saturated, we note
that by \cite[LG 4.26, Theorem 4]{Serre1}, the $p$-power map induces an isomorphism
$$
	H_\lambda\stackrel{\isom}{\longrightarrow} H_{\lambda+1}
$$
for all  $\lambda\in (\frac{1}{p-1},\infty)\cap v(\frm)$. As $H$ is complete, this implies that the
group is saturated.

The valuation on $R$ takes values in
$\frac{1}{e}\Z$ where $e$ is the ramification index. By definition
$H=H_\frac{\rho}{e}$.
We get
\[ \gr(H)\isom \bigoplus_{\lambda\in\frac{1}{e}\Na,\lambda\geq\frac{\rho}{e}} k^n\, . \]
 As an $\F_p[\epsilon]$-module $\gr(H)$ is freely generated by
an $\F_p$-basis of 
\[ \bigoplus_{\lambda\in\frac{1}{e}\Na,1+\frac{\rho}{e}>\lambda\geq\frac{\rho}{e}}k^n\]
This is finite because $[k:\F_p]<\infty$. 

If $e=1$, then only a single $\lambda$ occurs in the sum, namely $\rho$.
If $e>1$, then $1+\frac{\rho}{e}>\frac{\rho+1}{e}$ and the sum has generators in more than
one degree.
\end{proof}

An important example of the above construction arises from
separated smooth  group schemes $\bbG/R$.
The formal completion $\wh \bbG$ of $\bbG$ 
along its unit section is a formal group over $R$ 
and the associated standard
group $G$ is via $g\mapsto 1+g$ isomorphic to  
$$
	G\isom\ker (\bbG(R)\to \bbG(k)).
$$
\begin{example}\label{theexample}
To have an even more concrete example, consider $\bbG=\Gl_n$ over $R$.
Let $\pi$ be the uniformizer.
Then
$$
	\pi^\rho R=\left\{x\in R\mid v(x)>\frac{1}{p-1}\right\}.
$$
It follows from Lemma \ref{saturationlemma} that 
$$
	 H:=1+\pi^\rho M_n(R)\subset \Gl_n(R)
$$
is a saturated subgroup with respect to the filtration 
$\omega(1+(x_{i,j}))=\inf_{i,j}(v(x_{i,j}))$. As 
$$
\gr \left(  H\right)\isom\bigoplus_{\lambda\geq \frac{\rho}{e},\lambda\in\frac{1}{e}\Na} M_n(k)
$$
the  rank of $H$ is $n^2[R:\Z_p]$.
Note that this is not equi-$p$-valued for $e>1$.

However, we can view $\Gl_n(R)$ as 
the group of $\Z_p$-valued points of the Weil-restriction
$\bbG'=\Res_{R/\Z_p}\Gl_n$ which is 
a separated smooth group scheme over $\Z_p$
\cite[7.6, Proposition 5]{BLR}.
This point of view yields a different valuation on
the corresponding standard group, which is
\[ G'=\wh \bbG'(\Z_p)=1+pM_n(R)\subseteq\Gl_n(R).\]
By Lemma \ref{saturationlemma} $(G',\omega')$ is saturated and equi-$p$-valued
if $p>2$. 

As an explicit example, choose $R=\Z_p[\pi]$ with $\pi^2=p$ (hence $e=2$) and $n=1$.
Let $p>3$ for simplicity. 
$G=1+\pi R$ has rank $2$ with ordered basis
$x_1=1+\pi$, $x_2=1+\pi^2$ with
\[ \omega(x_1)=\frac{1}{2}, \omega(x_2)=1\]
On the other hand, $G'=1+pR$ has also rank $2$ with $x'_1:=1+p, x'_2:=1+\pi^3$
as an ordered basis. They satisfy
\[ \omega'(x'_1)=\omega'(x'_2)=1\]
$G'$ is saturated and equi-$p$-valued.
Compare this to
\[ \omega(x'_1)=1, \omega(x'_2)=\frac{3}{2}\]
under the inclusion $G'\subset G$.
\end{example}

\subsection{Valued rings, modules and the functor $\Sat$}
\begin{definition}\cite[I 2.1.1. and I 2.2.1]{L} 
A \emph{filtered} ring $\Omega$ is a ring together with a map
$$
	v:\Omega\to \bbR_+\cup \{\infty\}
$$
such that for $\lambda,\mu\in\Omega$
\begin{itemize}
	\item[1)] $v(\lambda-\mu)\ge \min(v(\lambda),v(\mu))$
	\item[2)]  $v(\lambda\mu)\geq v(\lambda)+v(\mu)$
	\item[3)]  $v(1)=0$.
\end{itemize}
Put $$
	\Omega_\nu:=\{\lambda\in\Omega\mid v(\lambda)\ge \nu\}.
$$
$\Omega$ is called \emph{valued} if  in addition
\begin{itemize}
	\item[2')]  $v(\lambda\mu) =  v(\lambda)+v(\mu)$
	\item[4)]The topology defined on $\Omega$ by the filtration
	$\Omega_\nu$ is separated.
\end{itemize}
\end{definition}
\begin{definition}\cite[I 2.1.3 and I 2.2.2]{L}
A \emph{filtered} module $M$ over a filtered ring $\Omega$ is an $\Omega$-module $M$ together with a map
$$
	w:M\to \bbR_+\cup \{\infty\}
$$
such that for $x,y\in M$ and $\lambda\in\Omega$
\begin{itemize}
	\item[1)] $w(x-y)\ge \min(w(x),w(y))$
	\item[2)]  $w(\lambda x)\geq v(\lambda)+w(x)$
\end{itemize}
Put $$
	M_\nu:=\{x\in M\mid w(x)\ge \nu\}.
$$
A filtered module over a valued ring $\Omega$ is called \emph{valued} if in addition
\begin{itemize}
	\item[2')] $w(\lambda x) =  v(\lambda)+w(x)$
	\item[3)]The topology defined on $M$ by the filtration
	$M_\nu$ is separated.
\end{itemize}
\end{definition}
Let $\Omega$ be a commutative valued ring and $A$ be an $\Omega$-algebra
(e.g. a Lie-algebra).
\begin{definition}\cite[I 2.2.4]{L}
An $\Omega$-algebra $A$ over a commutative valued ring $\Omega$ is \emph{valued}, if
it is valued as a ring and (with the same valuation map) 
valued as an $\Omega$-module.
\end{definition}
The following definition is an important technical tool
in Lazard's work.
\begin{definition}\cite[I 2.2.7]{L}
A valued module $M$ over a commutative valued ring $\Omega$ is called
\emph{divisible}, if for all $\lambda\in \Omega$ and $x\in M$ with $v(\lambda)\le w(x)$
there exists $y\in M$ such that $\lambda y=x$.
The module $M$ is \emph{saturated} if it is divisible and complete.
\end{definition}

Lazard shows in \cite[I 2.2.10]{L}  that the completion of a divisible module is
saturated. 
\begin{definition}\cite[I 2.2.11]{L}\label{satdefn}
The \emph{saturation} $\Sat M$ of a valued module $M$ over a commutative valued ring $\Omega$ 
is the completion of 
$$
	\div M:= \left\{y\in K\otimes_\Omega M\mid \tilde{w}(y)\ge 0\right\}.
$$
Here, $K$ is the fraction field of $\Omega$ and the valuation $w$ on $M$ is extended
to a map $\tilde{w}$ on 
$ K\otimes_\Omega M$ by $\tilde{w}(\lambda^{-1}\otimes m):= w(m)-v(\lambda)$ (which is well-defined,
see \cite[I 2.2.8]{L} ).
\end{definition}
The saturation $\Sat M$ satisfies the following universal property (\cite[I. 2.2.11.]{L}): For any
morphism $f:M\to N$ of $M$ into a saturated $\Omega$-module $N$ there is
a unique extension to a map $\Sat M\to N$.

\subsection{Group rings}
In this section we fix $\Omega=\Z_p$ with the standard valuation.
All algebras are over $\Z_p$.

For any group $G$ let $\bbZ_p[G]$ be the group ring with coefficients in $\bbZ_p$.
\begin{definition}\cite[II 2.2.1]{L}
Let $G$ be a pro-$p$-group. The {\em completed group ring} $\Z_p[[G]]$ is the projective limit 
$$
	\Z_p[[G]]:=\prolim \Z_p[G/U],
$$
where $U$ runs 
through all open normal subgroups of $G$ and every $\Z_p[G/U]$ carries the 
$p$-adic topology.
\end{definition}

In \cite{L} this ring is denoted $\mathrm{Al}G$.

\begin{definition}\cite[III 2.3.1.2]{L}\label{inducedfiltration}  Let $G$ be a $p$-filtered group. The
{\em induced filtration} $w$ on $\Z_p[G]$ is the lower bound for
all filtrations (as $\Z_p$-algebra) such that
\[ w(x-1)\geq \omega(x)\text{\ for all\ }x\in G\ .\]
\end{definition}

\begin{prop}\cite[III 2.3.3]{L}Let $G$ be $p$-valued. Then the induced filtration $w$ on $\Z_p[G]$ is a valuation (as $\Z_p$-module). If $G$ is compact
(or equivalently, pro-$p$), then $\Z_p[[G]]$ is the completion of $\Z_p[G]$
with respect to the valuation topology.
\end{prop}

\begin{example}\cite[V 2.2.1]{L}\label{exgroupring}
Let $G$ be $p$-valued, complete and of finite rank $d$. Let $\{x_i\}_{i=1,\dots,d}\subset G$ be an ordered basis of $G$. Then $\Z_p[[G]]$ admits 
$\{ z^\alpha\, |\,\alpha
\in\Na^d\}\subseteq\Z_p[[G]]$,
\[ z^\alpha:=\prod\limits_{i=1}^d (x_i-1)^{\alpha_i}\]
as a topological $\Z_p$-basis satisfying
\[ w(z^\alpha)=\sum\limits_{i=1}^d\alpha_i\omega(x_i).\]
The associated graded is $U_{\F_p[\epsilon]}\gr(G)$, the universal enveloping
algebra of the $\F_p[\epsilon]$-Lie-algebra $\gr(G)$.
\end{example}

\begin{remark}
Note that if $(G,\omega)$ is saturated and non-trivial,
then $\Z_p[[G]]$ is never saturated.
Indeed, since $\gr G\neq 0$ is a free $\F_p[\epsilon]$-module, we have $\gr^\nu G
\neq 0$ for arbitrarily large $\nu$, in particular there exists $g\in G$
with $\omega(g)\ge 1$. Then $x:=g-1\in\Z_p[[G]]$ satisfies $w(x)\ge 1=v(p)$, but
$x$ is not divisible by $p$ in $\Z_p[[G]]$.
\end{remark}

\begin{lemma}\label{completioniso}
The inclusion $\Z_p[G]\to \bbZ_p[[G]]$ induces an isomorphism
$$
\Sat \Z_p[G]\isom \Sat\bbZ_p[[G]].
$$
\end{lemma}
\begin{proof}
 By \cite[I 2.2.2]{L} the natural map $\Z_p[G]\to \Z_p[[G]]$ is injective.
It extends to 
\[ \Sat \Z_p[G]\to \Sat \Z_p[[G]] \]
 On the other hand,
$\Sat \Z_p[G]$ is complete, hence there is a natural map
$\Z_p[[G]]\to \Sat\Z_p[G]$. As the right hand side is saturated, it 
extends to
\[ \Sat \Z_p[[G]] \to \Sat \Z_p[G]\ .\]
 The two maps are inverse to each other.
\end{proof}

\subsection{Enveloping algebras}

Let $L$ be a valued $\Z_p$-Lie-algebra and $UL$ its
enveloping algebra over $\Z_p$.
\begin{defn}\cite[IV 2.2.1]{L} The {\em canonical filtration}
\[ w: UL\to \R_+\cup\{\infty\}\]
is the lowest bound of all filtrations on $UL$ turning it
into a valued $\Z_p$-algebra such that the canonical
map $L\to UL$ is a morphism of valued modules.
\end{defn}
\begin{lemma}\cite[IV 2.2.5]{L}
$UL$ equipped with the canonical filtration is a valued
$\Z_p$-algebra and the natural morphism
\[ U\gr(L)\to \gr(UL)\]
is an isomorphism.
\end{lemma}

\subsection{Group-like and Lie-algebra-like elements}

Everything in this section applies to $A=\Sat\Z_p[[G]]$ where 
$G$ is a $p$-valued pro-$p$-group. We fix $\Omega=\Z_p$ with its
standard valuation.

\begin{definition}\cite[IV 1.3.1]{L} \label{2.5.1}Let $A$ be a valued $\Z_p$-algebra with diagonal 
\[\Delta: A\to \Sat (A\tensor_{\Z_p} A)\]
(\cite[IV 1.2.3]{L})
and augmentation $\epsilon$. Then we define $\GG$, $\LL$, 
$\GG^*$ and $\LL^*$ by
\begin{align}
\GG &=\{ x\in A| \epsilon(x)=1, \Delta(x)=x\tensor x\} \\
\GG^*&=\{x\in \GG| w(x)> (p-1)^{-1}\}\\
\LL &= \{ x\in  A| \Delta(x)=x\tensor 1+1\tensor x\} \\
\LL^* &=\{x\in \LL| w(x)> (p-1)^{-1}\}
\end{align}
\end{definition}
These subsets have the following structures.
\begin{lemma}\cite[IV 1.3.2.1 and 1.3.2.2]{L}
 $\GG$ and $\GG^*$ are monoids with respect to the 
multiplication of $A$. If $A$ is complete, $\GG^*$
is a group and $\LL$ and $\LL^*$ are Lie-algebras. Moreover,
 $\LL=\div\LL^*$.
\end{lemma}
For saturated $\bbZ_p$-algebras $A$ we know much more:
\begin{theorem}\cite[IV 1.3.5]{L}\label{grouplikedescription}
Let $A$ be a saturated $\Z_p$-algebra with diagonal. 
\begin{enumerate}
\item
The exponential maps $\GG^*$ to $\LL^*$ and the logarithm maps $\LL^*$ to $\GG^*$. They
are inverse homeomorphisms.
\item
The Lie-algebra $\LL$ is saturated. It is the saturation of $\LL^*$.
\item
$\GG^*$ is a saturated group for the filtration $\omega(x)=w(x-1)$.
\item The associated graded $\gr\LL^*$ and $\gr\GG^*$ are canonically isomorphic via the logarithm map.
\item $\LL^*$ and $\GG^*$ generate the same saturated associative subalgebra of $A$. 
\end{enumerate}
\end{theorem}

This has the following consequence for the universal enveloping algebra $UL$ of a 
valued Lie-algebra $L$.
 
\begin{theorem}\cite[IV 3.1.2 and IV 3.1.3]{L}\label{universalsaturation} Let $L$ be a valued Lie-algebra over $\Z_p$
and $UL$ its universal enveloping algebra. Then
\setcounter{equation}{0}
\begin{align} 
\Sat UL&=\Sat U\Sat L\\
\LL \Sat UL&=\Sat L\\
\GG\Sat UL &=\GG^*\Sat UL
\end{align}
\end{theorem}

The next result concerns the saturation of the group ring $\bbZ_p[G]$ (or 
equivalently of $\bbZ_p[[G]]$ by Lemma \ref{completioniso}).

\begin{theorem}\cite[IV 3.2.5]{L} \label{mainisom}
Let $G$ be a saturated group and $A=\Sat\bbZ_p[G]$ then
$$
	\cG^*=G.
$$
Let $U\cL$ be the universal enveloping algebra of $\cL$, then
the canonical map
$$
	\Sat U\cL\to \Sat\bbZ_p[G]
$$
is an isomorphism.
\end{theorem}

We introduce
new terminology.

\begin{defn}\label{defintlaz}
Let $G$ be a saturated group. We call
\[ \LL^*(G)=\LL^*\subset \Sat\Z_p[G]\]
the \emph{integral Lazard Lie-algebra} of $G$.
\end{defn}
The last theorem then reads
\[ \Sat U\LL^*(G)\isom \Sat\Z_p[G].\]
\begin{example}\label{computelie}
Consider the saturated group $H=1+\pi^\rho M_n(R)$ from Example
\ref{theexample} and the algebra $\Sat\Z_p[H]$.
We claim that the Lie-algebra $\LL^*=\cL^*(H)$ is $\pi^\rho M_n(R)$
and $\cL=M_n(R)$. To see
this, note that by Theorem~\ref{mainisom} $H=\cG^*$ and that by 
Theorem \ref{grouplikedescription}
$\cL^*$ consists of the logarithms of $\cG^*$. By \cite[III 1.1.4 and 1.1.5.]{L}.
\begin{align} \Log: 1+\pi^\rho M_n(R)&\to \pi^\rho M_n(R)\ ; 1+x\mapsto \sum_{n\geq 1} (-1)^{n+1}\frac{x^n}{n}\\
\exp: \pi^\rho M_n(R)&\to 1+\pi^\rho M_n(R)\ ; x\mapsto \sum_{n\geq 0} \frac{x^n}{n!}
\end{align}
are both convergent and inverse to each other. By 
Theorem \ref{grouplikedescription}, $\cL$ is
the saturation of $\cL^*$, which is by Definition \ref{satdefn}
$$
	\cL=\{ x\in K\otimes_R \pi^\rho M_n(R)\mid w(x)\ge 0\}= M_n(R).
$$
\end{example}

\begin{example}\label{standardlie} In general the Lazard Lie-algebra does not coincide
with the algebraic Lie-algebra.
Let $\bbG$ be a separated smooth group scheme over $\Z_p$ and
$\Lie(\bbG)$ its $\Z_p$-Lie-algebra.
If $t_1,\dots,t_n$ are
formal coordinates of $\bbG$ around $e$, then $\frac{\partial}{\partial t_1},\dots,\frac{\partial}{\partial t_n}$
are a $\Z_p$-basis of $\Lie(\bbG)$.

Let $G$ be the associated standard group over $\Z_p$ as in Section \ref{standardsection}.
Let $H$ be the saturated subgroup of $G$, see Lemma \ref{saturationlemma}.
We have $H=G=(p\Z_p)^n$ if $p\neq 2$  and $H=(4\Z_2)^n$ if $p=2$.
Let $x_1,\dots,x_n$ be the standard ordered basis of $H$.
We put
\[ \delta_i=\log x_i\in\Sat\Z_p[[H]]. \]
By \cite[IV 3.3.6]{L} 
they form a $\Z_p$-basis of $\LL^*(H)$. 
As explained in \cite{HK} Section 4.2 and 4.3 the $\delta_i$ can be viewed as derivations on
$\Z_p[[t_1,\dots,t_n]]$, the coordinate ring of $\hat{\bbG}$.
Note, however, that the coordinate $\lambda_i$ in loc.cit.
takes values in all of $\Z_p$ on $G$. Hence $\lambda_i=pt_i$ for $p\neq 2$. This implies
\[ \delta_i=p\frac{\partial}{\partial t_i}\mid_{t=0}\]
Hence under the identification of \cite[Proposition 4.3.1]{HK}, we have
\[ \LL^*(H)=p\Lie(\bbG) . \]
For $p=2$ the argument gives
\[ \LL^*(H)=4\Lie(\bbG). \]
\end{example}

\subsection{Resolutions and Cohomology}

\begin{definition}\cite[I 2.1.16, 2.1.17]{L} 
Let $A$ be a filtered $\Z_p$-algebra, $M$ a filtered $A$-module. 
\begin{enumerate}
\item A family of $A$-linearly independent elements $(x_i)_{i\in I}$ of $M$ is called {\em filtered free}
if for every family $(\lambda_i)_{i\in I}$ of elements of $A$, almost all zero,
$$
	w\left(\sum_{i\in I}\lambda_ix_i\right)=\inf\limits_i (w(x_i)+v(\lambda_i)). 
$$
$M$ is called {\em filtered free} if it is generated by a filtered free family.
\item Suppose $A$ is complete. 
$M$ is called {\em complete free} if it is the completion of the submodule generated by a filtered free family.
\end{enumerate}
\end{definition}

If $A$ is complete and $M$ is filtered free of finite rank, then $M$ is also complete free. 

\begin{definition}\cite[V 1.1.3, 1.1.4, 1.1.7]{L}
Let $A$ be a filtered augmented $\Z_p$-algebra, $M$ a filtered $A$-module.
\begin{enumerate}
\item A {\em filtered acyclic resolution} $X_\bullet$ is a chain complex of filtered $A$-modules together with an augmentation $\epsilon:X_\bullet\rightarrow M$
such that for all $\nu\in \R_+$ the morphism
$\epsilon_\nu: X_{\bullet \nu}\to M_\nu$ is a quasi-isomorphism.
\item A {\em split filtered resolution} $X_\bullet$ of $M$ is a morphism $\epsilon:X_\bullet\to M$
of chain complexes of filtered $A$-modules together with filtered morphisms of $\Z_p$-modules (sic, not $A$-linear!)
\[ \eta:M\to X_0,\hspace{1cm} s_n:X_n\to X_{n+1}\]
defining a homotopy between $\mathrm{id}$ and $0$ on the extended complex $X_\bullet\stackrel{\epsilon}{\to} M$ and such that $s_0\eta=0$.
Note that a split filtered resolution is a filtered acyclic resolution.
\item Let $A$ be complete.
We call {\em complete free acyclic resolution} of 
$M$ a filtered acyclic resolution $X_\bullet$ by complete free modules.
\item Let $X_\bullet$ be a complete free acyclic resolution of the trivial $A$-module $\Z_p$, $M$ a complete
$A$-module with linear topology. We call
\[ H^n_c(A,M)= H^n(\Hom_c(X_\bullet,M))\]
(with $\Hom_c$ continuous $A$-linear maps) the {\em $n$-th continuous cohomology} of $A$ with coefficients in $M$.
\item Let $A$ be an augmented $\bbZ_p$-algebra, $M$ an $A$-module. We call
\[ H^n(A,M)=\Ext_A^n(\bbZ_p,M)\]
the {\em $n$-th cohomology} of $A$ with coefficients in $M$.
\end{enumerate}
\end{definition}


\section{An integral version of the Lazard isomorphism}\label{intsection}

The purpose of this section is to establish that continuous group cohomology and Lie-algebra cohomology agree
with integral coefficients, at least under certain technical assumptions.
This generalizes Lazard's result for coefficients in $\bbQ_p$-vector spaces.

\subsection{Results}

We fix a saturated and compact group $(G,\omega)$ of finite rank $d$. In particular,
$G$ is a pro-$p$-group by Proposition \ref{pro-p-crit}. We assume
\begin{itemize}
\item $(G,\omega)$ is equi-$p$-valued
\item $\omega$ takes values in $\frac{1}{e}\Z$.
\end{itemize}
Recall that the integral Lazard Lie-algebra 
\[ \LL^*(G)=\cL^*\Sat\Z_p[[G]],\]
is a finite free $\Z_p$-Lie-algebra.

For technical reasons we fix a totally ramified extension $\Q_p\subseteq K$ of degree $e$ with ring
of integers $\Oh\subseteq K$, uniformizer $\pi\in\Oh$. The valuation on
$\Oh$ is normalized by $v(p)=1$. 

Let $M$ be a linearly topologized complete $\Z_p$-module with a continuous, $\Z_p$-linear action of $G$.
Thus, $M$ is a $\Z_p[[G]]$-module (\cite[II 2.2.6]{L}). We assume that
\begin{itemize}
\item the  $\Z_p[[G]]$-module structure on $M$  extends to a $\Sat \Z_p[[G]]$-module structure.
\end{itemize}
$M$ is canonically a $\LL^*(G)$-module.

We are going to prove in Section \ref{sectionproof}:
\begin{thm}\label{intlazard}
Let $(G,\omega)$ and $M$ be as above, then:
\begin{enumerate}
\item
There is an isomorphism of graded $\Oh$-modules
\[ \phi_G(M):H^*_c(G,M)\otimes_{\bbZ_p} \Oh\simeq H^*(\LL^*(G),M)\otimes_{\bbZ_p}\Oh.\]
It is natural in $M$.
\item
If in addition $M$ is a $\Q_p$-vector
space, then this isomorphism agrees with Lazard's in \cite[V 2.4.9.]{L} 
\item Let $H$ be another group satisfying the assumptions of the Theorem and
$f:G\to H$ a group homomorphism filtered for the chosen filtrations. In addition
assume that $\gr(H)$ is generated in degree $\frac{1}{e}$. Then the isomorphism
is natural with respect to $f$.
\item
If $\gr(G)$ has generators in degree $\frac{1}{e}$, then
the isomorphism is compatible with cup-products as follows:\\
Assume that $M',M''$ satisfy the same assumptions as $M$ does and that
\[ \alpha: M\hat{\otimes}_{\Z_p} M'\to M''\]
is $\Sat\Z_p[[G]]$-linear. Then the diagram
\[ 
\xymatrix{ \left(H^*_c(G,M)\otimes_{\Z_p} H^*_c(G,M')\right) \otimes\Oh \ar[r] \ar[d]^{\phi_G(M)\otimes\phi_G(M')} & H^*_c(G,M'')\otimes\Oh\ar[d]^{\phi_G(M'')}\\
\left(H^*(\LL^*(G),M)\otimes_{\Z_p} H^*(\LL^*(G),M')\right)\otimes\Oh \ar[r] & H^*_c(\LL^*(G),M'')\otimes \Oh}
\]
commutes. Here, the horizontal maps are the $\Oh$-linear extensions of the cup-product defined by $\alpha$.
\end{enumerate}
\end{thm}

\begin{remark}
\begin{itemize}
\item[i)] If $H^*(\LL^*(G),M)$ is a finitely
generated $\Z_p$-module, e.g. when $M$ is of finite type,
then this implies by 
the structure of finitely generated modules over principal ideal domains
the existence of an isomorphism of graded $\Z_p$-modules
\[ H^*_c(G,M)\simeq H^*(\LL^*(G),M).\]
However, it is not clear if this isomorphism is natural or
compatible with cup-products.
\item[ii)] According to \cite[V 2.2.6.3 and 2.2.7.2]{L} the mod $p$ cohomology
of an equi-$p$-valued group $G$ is simply an exterior algebra
\[ H^*_c(G,\F_p)=\Lambda^*_{\F_p}(H^1_c(G,\F_p))\]
but the cohomology with torsion free coefficients is more interesting, e.g.
if $G$ is not abelian then the $\Q_p$-Betti numbers of $G$ are different
from the $\F_p$-Betti numbers showing
that $H^*_c(G,\Z_p)$ contains non-trivial torsion.
\end{itemize}
 
\end{remark}

It is not  obvious to see which groups satisfy the assumptions
of Theorem~\ref{intlazard}. We discuss in section \ref{examplesection} standard groups and 
uniform pro-$p$-groups, which satisfy the assumptions of Theorem \ref{intlazard}.
The next section~\ref{explanations} discusses  the assumptions with
some examples.
\subsection{Some examples concerning the assumptions of Theorem \ref{intlazard}}\label{explanations}
Here we illustrate the assumptions of Theorem \ref{intlazard}
by a series of remarks and examples.

The integral Lazard isomorphism may not hold for all
topologically finitely generated pro-$p$-groups without $p$-torsion. 
However, the assumptions of the Theorem are too restrictive.

\begin{example} \label{morava}
Assume $p\ge 5$ and let $D/\Q_p$ be the quaternion-algebra,
$\Oh\subseteq D$ its maximal order and $\Pi\in \Oh$ a prime element.
Using \cite[II 1.1.9]{L} one can check that \[ G:=1+\Pi\Oh\subseteq\Oh^*\]
is $p$-saturated. From \cite[Theorem 6.3.22]{R} or \cite[Proposition 7]{Henn}, we know that 
\[ \dim_{\F_p}\, H^i_c(G,\F_p)=1,3,4,3,1\,\, (0\leq i\leq 4).\]
In particular, $H^*_c(G,\F_p)\neq\Lambda^* H^1_c(G,\F_p)$ and $G$
does not admit an equi-$p$-valuation by \cite[V 2.2.6.3 and 2.2.7.2.]{L} 
However, one can
by direct arguments establish an isomorphism
\[ H^*_c(G,\F_p)\simeq H^*(\LL^*(G),\F_p)\]
of graded $\F_p$-algebras; cf. Remark \ref{superfluous}. The proof of \cite[Proposition 7]{Henn} shows
that the same result holds for coefficients in $\bbZ_p$.
\end{example}

Not even saturatedness is necessary.
\begin{example}
Let $G=1+p^2\Z_p$ for $p\neq 2$. This group is not saturated 
for the obvious filtration, rather we have
\[\Sat(G)=1+p\Z_p. \]
Put 
\[ \LL^*(G)\subset \LL^*(\Sat(G))\]
the image of $G$ under the logarithm map. We still get an isomorphism
\[ H^*(G,\Z_p)\to H^*(\LL^*(G),\Z_p)\]
induced by $\log$. It is compatible with the one for $\Sat(G)$. 
\end{example}

\begin{rem}\label{applicability}
\begin{enumerate}
\item We are unaware of a group-theoretical characterization of those 
pro-$p$-groups satisfying the assumption of Theorem \ref{intlazard}, but
the remark on page 163 of \cite{ST} suggests that they are 
closely related to uniform pro-$p$-groups. 
\item It is in general difficult to decide if a given $\Z_p[[G]]$-module
structure extends over $\Sat\Z_p[[G]]$, and we refer to \cite[page 200]{T}
and especially to the proof of \cite[Corollary 9.3]{T} for further discussion 
and useful sufficient conditions.
\item In Theorem \ref{uniform} we establish a sufficient
condition for both problems which have to be addressed here.
\end{enumerate}
\end{rem}
There are examples of groups which are saturated with respect to one filtration
but not with respect to another. It can also happen that the group is saturated with respect to two filtrations but
only equi-$p$-valued for one of them.

\begin{example} Let $K/\Q_p$ be a finite extension
with ramification index $e$. Let $\Oh$ be its
ring of integers with uniformizer $\pi$. 
As discussed in Example \ref{theexample} the group
\[ 1+p M_n(\Oh)\]
carries two natural filtrations $\omega$ and $\omega'$.
Recall that $\rho$ is the smallest integer bigger than $\frac{e}{p-1}$.
\begin{enumerate}
\item If $p=5$, $e=2$, then $\rho=1$ 
and hence $\pi^\rho\neq 5$. This implies that $1+5M_n(\Oh)$
is saturated with respect to $\omega'$ but not with respect to
$\omega$.
\item If $p=3$, $e=2$, then $\rho=2$ and hence $\pi^\rho=3$.
The group $1+3M_n(\Oh)$ is saturated with respect to $\omega$ and
$\omega'$, but but only equi-$3$-valued with respect
to the second.
\end{enumerate}
\end{example}

\subsection{The case of standard groups and uniform pro-$p$-groups}\label{examplesection}
We discuss two examples, where the assumptions of Theorem \ref{intlazard} are satisfied.
First we consider standard groups and then uniform pro-$p$-groups.
\begin{example}
Let $\G/\Z_p$ be a separated smooth group scheme, $G=\ker(\G(\Z_p)\to \G(\F_p))$ the associated
standard group (see Section \ref{standardsection}).
Its filtration takes values in $\Z$.
By Lemma \ref{saturationlemma}, there is an open subgroup $H$
of $G$ which is saturated and equi-$p$-valued. 
If $p\neq 2$, then $H=G$ and the generators have degree $1$.
If $p=2$, then the generators have degree $2$.
$H$ satisfies the assumptions of the Theorem with $e=1$.

Let $M=\Z_p$ with the trivial operation of $H$. It also satisfies
the assumptions of the Theorem. Hence there is a natural
isomorphism of graded $\Z_p$-modules
\[ H^*_c(H,\Z_p)\simeq H^*(\LL^*(H),\Z_p).\]
For $p\neq 2$ it is compatible with cup-products.
\end{example}
This example generalizes to a larger class of groups.
Recall the notion of a {\em uniform} or {\em uniformly
powerful} pro-$p$-group 
from \cite[Definition 3.1 and Definition 4.1]{DDMS}:
\begin{definition} \label{powerful}
A pro-$p$-group $G$ is {\em uniform}
if \begin{itemize}
\item[i)] $G$ is topologically finitely generated.
\item[ii)] For $p\neq 2$ (resp. $p=2$), $G/\overline{G^p}$ (resp. $G/\overline{G^4}$) is abelian.
\item[iii)] Denoting $G=G_1\supseteq G_2\supseteq\ldots$ the lower $p$-series
of $G$, we have $[G_i:G_{i+1}]=[G_1:G_2]$ for all $i\ge 2$.
\end{itemize}
\end{definition}
To understand what is special about $p=2$ here, note that the pro-$2$-group
$\Z_2^*=1+2\Z_2$ is not uniform but $1+4\Z_2$ is.\\
 
The Lie-algebra $\g$ of a uniform pro-$p$-group $G$
is constructed in \cite[\S 8.2]{DDMS} and coincides
with the integral Lazard Lie algebra $\LL^*(G)$ by \cite[Lemma 8.14]{DDMS}.
\begin{thm}\label{uniform}
Let $p\neq 2$ (resp. $p=2$) be a prime, $G$ a uniform pro-$p$-group and $M$
a finite free $\Z_p$-module with a continuous action of $G$ such that the
resulting group homomorphism
\[ \varrho:G\longrightarrow\Aut_{\Z_p}(M)\]
has image in $1+p\End_{\Z_p}(M)$ (resp. in $1+4\End_{Z_2}(M)$).
Then $M$ is canonically a module for the Lie-algebra $\g$ of $G$
and there is an isomorphism of graded $\Z_p$-modules 
\begin{equation}\label{ersteriso}
 H^*_c(G,M)\simeq H^*(\g,M)
\end{equation}
which in case $p\neq 2$ 
is compatible with cup-products whenever these are defined.
\end{thm}

\begin{rem} If $G$ is an arbitrary $\Q_p$-analytic group acting
continuously on the finite free $\Z_p$-module $M$, then there are 
arbitrarily small open subgroups $U\subseteq G$ such that the action of $U$ on $M$
satisfies the assumptions of Theorem \ref{uniform}.
\end{rem}

\begin{proof} We have the following two claims for $p\neq 2$ (resp. $p=2$):
\begin{enumerate}
\item  $G$ admits a valuation $\omega$ for which it is $p$-saturated of finite rank
and equi-$p$-valued with an ordered basis consisting of elements of 
filtration $1$ (resp. filtration $2$).\\
\item The $\Z_p$-module $M$ admit a valuation $w$ for which it is saturated and
such that for all $g\in G,m\in M:w((g-1)m)\ge w(m)+\omega(g)$.
\end{enumerate}
Granting these claims, we see as in \cite[pages 200-201]{T} that
the $\Z_p[[G]]$-module structure of $M$ extends over $\Sat\Z_p[[G]]$ and 
hence obtain (\ref{ersteriso}) by applying Theorem \ref{intlazard},(1) with $\Oh=\Z_p$
and remarking that $\g=\LL^*((G,\omega))$.\\
The proof of claim (1) is essentially given in \cite[Remark on page 163]{ST}
but we include details for convenience:
The lower $p$-series \[ G=G_1\supseteq G_2\supseteq\ldots\]
\cite[Definition 4.1]{DDMS} consists of normal subgroups satisfying
$ (G_n,G_m) \subseteq G_{n+m}$ and $\cap_{n\ge 1} G_n=\{ e \}$
\cite[Proposition 1.16]{DDMS} and hence
\[ \omega(x):=\sup\{ n\in\Na\, | \, x\in G_n\}, x\in G\]
defines a filtration of $G$ by \cite[II.1.1.2.4.]{L}.
Now \cite[Lemma 4.10]{DDMS} states that for all $n,k\geq 1$ the 
$p^n$-th power map of $G$ is a homeomorphism $G_k\stackrel{\isom}{\to}G_{k+n}$
and induces bijections $G_k/G_{k+l}\stackrel{\isom}{\to} G_{k+n}/G_{k+n+l}$
for all $l\ge 0$. 

In case $p\neq 2$, we get from this all the properties of $\omega$ in Definition \ref{filtrationdefn}
we need:
\begin{enumerate}
\item[4)] is trivial since $1>\frac{1}{p-1}$.
\item[5)] If $x\in G$
has filtration $n=\omega(x)$ then $[x^p]\in G_{n+1}/G_{n+2}$ is non-trivial,
i.e. $\omega(x^p)=\omega(x)+1$.
\item[6)]
If $x\in G$ satisfies 
$\omega(x)>1+\frac{1}{p-1}$ then $\omega(x)\ge 2$, hence $x\in G^p$.
\end{enumerate}
As $G$ is complete, we see that $(G,\omega)$
is $p$-saturated, clearly of finite rank. More precisely, from the above we
get that $\gr G$ is $\F_p[\epsilon]$-free on $\gr^1 G$ and thus $G$ is 
equi-$p$-valued with an ordered basis consisting of elements of filtration $1$.
This settles claim (1) in case $p\neq 2$.\\
In case $p=2$, $\omega$ satisfies all conditions in Definition
\ref{filtrationdefn}, except 4), so $(G,\omega)$ is in particular $2$-
filtered and $\gr G$ has a structure of mixed Lie algebra over $\F_2$ 
\cite[II 1.2.5]{L}. Note that the only $\nu\in\R^+$ with $\nu\leq\frac{1}{p-1}=1$
and $\gr_\nu G\neq 0$ is $\gr_1 G$. From Definition \ref{powerful}, ii)
we have $[\gr_1 G,\gr_1 G]=0$ which easily implies that $\gr G$ is abelian.
Since $\omega$ has integer values, this means that 

\[ \omega([x,y])\ge\omega(x)+\omega(y)+1\, ,\, x,y\in G.\]

Using this, it is easy to see that $\omega':=\omega+1$ is a filtration of $G$
with the properties stated in claim (1) in case $p=2$.\\

As for claim (2), $p$ being arbitrary now,
we choose a $\Z_p$-basis $\{e_i\}\subseteq M$ and declare
it to be a filtered basis with $w(e_i)=0$, i.e. 
\[ w\left(\sum_i\lambda_i e_i\right)=
\inf_i \{v(\lambda_i)\},\  \text{for\ } \lambda_i\in \Z_p\]
Clearly, $(M,w)$ is saturated. Assume $p\neq 2$.

We consider the continuous homomorphism of pro-$p$-groups
\[ 
	\varrho: G\longrightarrow 1+p\End_{\Z_p}(M)=:G'
\]
and claim that the lower-$p$-series of $G'$ is given by 
$G'_n=1+p^n\End_{\Z_p}(M)$, $n\ge 1$. Since $G'$ is powerful, \cite[Lemma 2.4]{DDMS}
gives $G'_{n+1}=\Phi(G'_n)$, the Frattini subgroup, for all $n\ge 1$ and
arguing inductively, it suffices to see that $\Phi(1+p^n\End_{\Z_p}(M))=
1+p^{n+1}\End_{\Z_p}(M)$. Since the Frattini subgroup is generated by $p$-th
powers and commutators, we have ``$\supseteq$'' and \cite[Proposition 1.16]{DDMS}
then gives
\begin{align*}
 \Phi(1+p^n\End_{\Z_p}(M))&/(1+p^{n+1}\End_{\Z_p}(M))\\
 &=\Phi((1+p^n\End_{\Z_p}(M)))/(1+p^{n+1}\End_{\Z_p}(M))\\
 &=\Phi((\F_p,+)^{n^2})=0.
\end{align*}
Since $\varrho$ respects the lower $p$-series, we conclude that 
\[ 
\varrho(G_n)\subseteq 1+p^n\End_{\Z_p}(M)\, , \, n\ge 1
\]
which implies that $w((g-1)m)\ge w(m)+\omega(g)$ for all $g\in G,m\in M$.\\
This settles claim (2) in case $p\neq 2$ and the argument in case $p=2$ is 
an obvious modification which we leave to the reader.\\

Finally, to see compatibility with cup products,
assume $p\neq 2$ and $M',M''$ satisfy the same assumptions as $M$ does and 
\[ \alpha:M\otimes_{\Z_p} M'\to M''\]
is $G$-linear defining cup-products in $H^*_c(G,-)$. Then both the source
and the target of $\alpha$ are canonically $\Sat\Z_p[[G]]$-modules
as seen above and $\alpha$ is $\Sat\Z_p[[G]]$-linear. Hence (\ref{ersteriso})
is compatible with cup-products by Theorem \ref{intlazard},(4).
\end{proof}

\subsection{Proof of Theorem \ref{intlazard}}\label{sectionproof}
We now describe the set-up for the rest of the section.

We fix a saturated group $(G,\omega)$ of finite rank $d$. Let
\[ \LL^*(G)=\cL^*\Sat\Z_p[[G]] \]
be its integral Lazard Lie-algebra. It is a finite free $\Z_p$-module.

We fix an ordered basis $\{x_1,\ldots,x_d\}\subseteq G$, and put 
$\omega_i:=\omega(x_i)$.
For every $0\leq k\leq n$ let
\[ \II_k:=\{ (i_1,\ldots,i_k)\, |\,1\leq i_1<\ldots<i_k\leq n\}\]
and for $I\in\II_k$ write $|I|:=\sum\limits_{s=1}^k\omega_s$.
For $I\in \II_0=\emptyset$ we put by abuse of notation $|I|=0$.

We assume that there
exists an integer $e\ge 1$ such that $\omega(G)\subseteq\frac{1}{e}\Z$
and fix a totally ramified extension $\Q_p\subseteq K$ of degree $e$ with ring
of integers $\Oh\subseteq K$, uniformizer $\pi\in\Oh$. The valuation on
$\Oh$ is normalized by $v(p)=1$. The artificial introduction of $\cO$ is a trick 
invented by Totaro in \cite{T}.
In this section all valued modules and algebras
are over $\Oh$. In particular, the saturation functor is taken in the
category of valued $\Oh$-modules. 

The inclusion $\Z_p\subseteq \Oh$ induces
\[
\F_p[\epsilon]=\gr \Z_p\subseteq\gr\Oh=\F_p[\epsilon_K]
\]
where $\epsilon$ (resp. $\epsilon_K$) is the leading term of
$p\in\Z_p$ (resp. $\pi\in\Oh$). We have $\epsilon_K^e\in\F_p^*\cdot\epsilon$,
in particular the degree of $\epsilon_K$ is $\frac{1}{e}$.

If $M$ is a valued $\Oh$-module, $\gr(M)$ is canonically an $\F_p[\epsilon_K]$-module.
As pointed out by Totaro (\cite[p. 201]{T}) it follows directly from the
definitions that
\[
\gr(\Sat(M))=\left(\gr M\otimes_{\F_p[\epsilon_K]}\F_p[\epsilon_K^{\pm 1}]\right)_{\mbox{\tiny degree}\geq 0}\ .
\]
Let
\[ A:=\Oh[[G]]:=\prolim\limits_{U\subseteq G\mbox{\tiny\ open normal}} \Oh[G/U],\]
and
\[ B:=U_\Oh(\LL^*(G)\otimes_{\Z_p}\Oh)^\wedge=U_{\Z_p}(\LL^*(G))^\wedge\otimes_{\Z_p}\Oh,\]
the completion of the universal enveloping algebra with respect to 
its canonical filtration. (This filtration is easily seen (using
Poincar\'e-Birkhoff-Witt) to be the $p$-adic filtration, hence the claimed
equality because $\Oh$ is finite free as a $\Z_p$-module.) We finally introduce
using Theorem \ref{mainisom}
\[ C:=\Sat A\isom \Sat B.\]

\begin{lemma} $\gr(A)=\gr(B)$ inside $\gr(C)$.\end{lemma}
\begin{proof}
We have
\[
	\gr A=\gr(\Z_p[[G]]\otimes_{\Z_p}\Oh)=
	U_{\F_p[\epsilon]}(\gr G)\otimes_{\F_p[\epsilon]}\F_p[\epsilon_K]=
	U_{\F_p[\epsilon_K]}(\gr G\otimes_{\F_p[\epsilon]}\F_p[\epsilon_K]),
\]
whereas 
\[
\gr B=\gr(U_{\Z_p}(\LL^*(G))\otimes_{\Z_p}\Oh)=U_{\F_p[\epsilon_K]}(\gr \LL^*(G)\otimes_{\F_p[\epsilon]}\F_p[\epsilon_K]).
\]
Since $\gr \LL^*(G)=\gr G$ by Theorem \ref{mainisom} and 
Theorem \ref{grouplikedescription}, (4) the claim follows.
\end{proof}
\begin{rem} Moreover, Totaro shows 
in \cite[pages 201-202]{T}
that for 
\[ \mathfrak{t}:=(\gr G\otimes_{\F_p[\epsilon]}\F_p[\epsilon_K^{\pm 1}])_{\mbox{\tiny degree}\ge 0},\]
a finite graded free $\F_p[\epsilon_K]$-Lie-algebra with generators in 
degree zero, we have $\gr C=U_{\F_p[\epsilon_K]}(\mathfrak{t})$.
\end{rem}

\begin{lemma}\label{modulesat}
\begin{enumerate}
\item Let $X$ be a filtered free $A$-module  with $A$-basis 
$e_1,\dots,e_r$. Then $\Sat X$ is a filtered free $C$-module on
generators 
\[ e'_i=\pi^{-e w(e_i)}e_i,\ i=1,\dots,r\]
\item Let $Y$ be a filtered free $B$-module with $B$-basis
$f_1,\dots,f_s$. Then $\Sat Y$ is a filtered free $C$-module on
generators 
\[f'_j=\pi^{-e w(f_j)}f_j,\ j=1,\dots,s\]
\end{enumerate}
\end{lemma}
\begin{proof}
It suffices to consider the case of the algebra $A$. The argument for
$B$ is the same. Without loss of generality $r=1$. 
By construction (and because $X$ is torsion-free), there are embeddings 
\[\xymatrix{
&\div X\ar[rd]\\
X\ar[ru]\ar[rd]&& K\otimes_\Oh X\\
&(\div A)\otimes_AX\ar[ru]
}\]
By assumption, any element $x$ of $K\otimes X$ can be written in the form
\[ x= \pi^v a e_1=\pi^{v+e w(e_1)}ae'_1, \in\Z, a\in A\]
It is in $\div X$ if and only if 
\[ w(x)=\frac{v}{e}+w(a)+w(e_1)\geq 0
\]
This is equivalent to $\pi^{v+ew(e_1)}a\in\div A$ and hence $x\in\div(A)e'_1$.
Hence $\div X=(\div A)\otimes_A X$.

Finally, apply the completion functor to the equality.
\end{proof}
\begin{rem}
This is the step where we make use of the coefficient extension
to $\Oh$.
\end{rem}

Both $A$ and $B$ are canonically subrings of $C$, and our first aim
is to compare the cohomology of the (abstract) rings $A$ and $B$
with that of $C$. Both $A$ and $B$ are augmented $\Oh$-algebras, hence we have an $A$- (resp. $B$-)
module structure on $\Oh$ which we will refer to as trivial.

\begin{prop}\label{resolution}
\begin{enumerate}
\item The trivial $A$-module $\Oh$ admits a resolution
$X_\bullet$ such that $X_k$ is filtered free of rank ${d\choose k}$
over $A$ on generators $\{e_I\, |\, I\in\II_k\}$ of filtration $w(e_I)=|I|$.
\item The trivial $B$-module $\Oh$ admits a resolution
$Y_\bullet$ such that $Y_k$ is filtered free of rank ${d\choose k}$
over $B$ on generators $\{f_I\, |\, I\in\II_k\}$ of filtration $w(f_I)=|I|$.
\item
Furthermore, $X_\bullet$ and $Y_\bullet$ can be chosen such that 
$\gr X_\bullet=\gr Y_\bullet$ as complexes of $\gr A=\gr B$-modules.
\end{enumerate}
\end{prop}

\begin{proof}
\begin{enumerate}\item The base extension from $\Z_p$
to $\Oh$ of the quasi-minimal complex of $G$ has the desired properties \cite[V 2.2.2.]{L}.  To see that the generators have the indicated
filtration, remember that the quasi-minimal complex is obtained by lifting the standard complex $\overline{X}_\bullet$ of the $\F_p[\epsilon]$-Lie-algebra $\gr(G)$ which has $\overline{X}_k=\Lambda^k_{\F_p[\epsilon]}(\gr(G))$ finite graded free on 
\[
\{ x_{i_1}G_{\omega_1}^+\wedge\ldots\wedge x_{i_k}G_{\omega_k}^+\}.
\]
\item The Lie-algebra $\LL^*(G)$ is $\Z_p$-free on generators $\Log(x_i)$ of filtration $\omega_i$. Hence the standard complex of $\LL^*(G)\otimes_{\Z_p}\Oh$ is as desired.
\item The equality $\gr X_\bullet=\gr Y_\bullet$ follows by construction
from $\gr(G)\isom \gr\LL^*(G)$.
\end{enumerate}
\end{proof}

\begin{example}
If $G$ is equi-$p$-valued, i.e, $\omega_i=\omega_j$ for all $i,j$, 
then $X_\bullet$ and $Y_\bullet$ are {\em minimal} in the sense of
\cite[V 2.2.5]{L}, i.e., $X_\bullet\otimes\F_p$ and $Y_\bullet\otimes \F_p$ have zero differentials.
\end{example}

For the following, we fix complexes $X_\bullet$ and $Y_\bullet$
satisfying the conclusion of Proposition \ref{resolution}.
Note that $C$ is an augmented $\Oh$-algebra with augmentation extending 
both the one of $A$ and the one of $B$.

\begin{lemma}\label{structuresat} Both $\Sat X_\bullet$ and $\Sat Y_\bullet$
are finite  filtered  resolutions of the trivial $C$-module $\Oh$, the modules
$Sat X_k$ (resp. $\Sat Y_k$) being filtered free on generators
$\{ \pi^{-e|I|}e_I\, |\, I\in\II_k\}$ (resp. $\{ \pi^{-e|I|}f_I\, |\, I\in\II_k\}$) of filtration zero over $C$. 
\end{lemma}

\begin{proof} Clearly, $\Sat X_\bullet$ and $\Sat Y_\bullet$
are canonically complexes of $C$-modules. Since both $X_\bullet$
and $Y_\bullet$ admit the structure of a {\em split} resolution, and
this structure is preserved by the additive functor $\Sat$, 
both $\Sat X_\bullet$ and $\Sat Y_\bullet$ are resolutions of $\Sat\Oh=\Oh$.

The statement on generators follows directly from Lemma \ref{modulesat}.
\end{proof}

For $0\leq k\leq n$ and $I\in\II_k$ denote by $e'_I\in\Sat X_k$
(resp. $f'_I\in \Sat Y_k$) the $C$-generators found above, i.e.
$e'_I:=\pi^{-e|I|}e_I, f'_I:=\pi^{-e|I|}f_I$.\\
We see that the canonical morphisms of complexes over $C$
\[ C\otimes_A X_\bullet\hookrightarrow\Sat X_\bullet\]
and
\[ C\otimes_B Y_\bullet\hookrightarrow\Sat Y_\bullet\]
are injective.

We pause to remark that, evidently, the above injections are isomorphisms rationally, a key input in Lazard's comparison isomorphism
for rational coefficients.
Similarly, an integral version of this
comparison isomorphism is essentially equivalent to $C\otimes_A X_\bullet$
being isomorphic to $C\otimes_B Y_\bullet$ and we proceed to prove this
in a special case as follows.

\begin{prop}\label{iso} There exists an isomorphism
 \[ \phi:\Sat X_\bullet\to\Sat Y_\bullet\] of filtered complexes 
over $C$ such that $\gr \phi=\id$ and $H^0(\phi)$ is the identity of $\Oh$.
Any two such $\phi$ are chain homotopic where the homotopy $h$ can be chosen such
that $\gr(h)=0$.
\end{prop}

\begin{proof} In order to construct the isomorphism it suffices
using \cite[V 2.1.5]{L}  (applicable by Proposition \ref{structuresat})
to canonically identify the complexes $\gr\Sat X_\bullet$ and
$\gr\Sat Y_\bullet$ of $\gr C$-modules. Recall from Proposition
\ref{resolution} that
\[ \gr X_\bullet=\gr Y_\bullet\]
This implies
\begin{multline*}
\gr\Sat X_\bullet=(\gr X_\bullet\otimes_{\F_p[\epsilon_K]}\F_p[\epsilon_K^{\pm 1}])_{\mbox{\tiny degree}\ge 0}\\
=(\gr Y_\bullet\otimes_{\F_p[\epsilon_K]}\F_p[\epsilon_K^{\pm 1}])_{\mbox{\tiny degree}\ge 0}=\gr\Sat Y_\bullet.
\end{multline*}

Now, $H^0(\phi)$ is an $\Oh$-linear automorphism of $\Oh$, hence given by multiplication with a unit $\alpha\in\Oh^*$. Using that its associated graded is the 
identity, one easily obtains $\alpha=1$, as claimed.

We turn to the construction of the homotopy. Let $\phi$, $\phi'$ be isomorphisms as above.
Let $e'_I\in \Sat(X_0)$ be a basis element. We need
to define $h_0(e_I)\in\Sat(Y_0)$ such that
\[ d h_0(e'_I)=(\phi-\phi')(e'_I)=:y_I\]
By assumption $\gr(\phi-\phi')=0$, and hence $y_I\in \Sat(Y_0)_{\frac{1}{e}}$.
As $\phi$ and $\phi'$ are isomorphisms of resolution of $\Oh$, we have
$\epsilon(y_I)=0$. Recall that $\Sat X_\bullet$ and $\Sat Y_\bullet$ are {\em filtered}
resolutions. Hence $y_I$ has a preimage $\tilde{y}_I\in\Sat(Y_1)_{\frac{1}{e}}$. Put
\[ h_0(e'_I)=\tilde{y}_I\]
By $C$-linearity, this defines $h_0$. It satisfies $\gr(h_0)=0$.
As usual, the same argument can be used inductively to
define $h_i$ for all $i\geq 0$.
\end{proof}

\begin{prop}\label{restrictiso} If, in the situation of 
Proposition \ref{iso}, $(G,\omega)$ is assumed to be equi-$p$-valued, 
then $\phi$ restricts to an isomorphism
\[ \psi:C\otimes_AX_\bullet\to C\otimes_BY_\bullet\]
of complexes over $C$.
If moreover $\gr(G)$ is generated in degree $\frac{1}{e}$, then
any two such isomorphisms are homotopic.
\end{prop}
\begin{proof}
We have the solid diagram of complexes over $C$
\[ 
\xymatrix{ C\otimes_A X_\bullet \ar@{^(->}[r]^{\iota_1} \ar@{.>}[d]^\psi & \Sat X_\bullet\ar[d]^\phi_{\simeq} \\
C\otimes_B Y_\bullet \ar@{^(->}[r]^{\iota_2} & \Sat Y_\bullet } 
\]
Since the horizontal maps are injective, $\phi$ factors as a chain-map
if for every $0\leq k\leq n$ we have 
\begin{equation}\tag{*} 
	\phi_k(C\otimes_A X_k)\subseteq C\otimes_B Y_k.
	\end{equation}

If $\psi$ exists, it is necessarily an isomorphism by completeness and
the fact that its associated graded is the identity. 
Alternatively observe that the following argument applies likewise to
$\phi^{-1}$ to produce an inverse of $\psi$.

To see what $(*)$
means, fix $0\leq k\leq n$ and remember the $C$-generators
$e_I\in C\otimes_A X_k,e'_I\in\Sat X_k, f_I\in C\otimes_B Y_k$
and $f'_I\in\Sat Y_k$ ($I\in\II_k$) satisfying $\iota_1(e_I)=\pi^{e|I|}e'_I$ and $\iota_2(f_I)=\pi^{e|I|}f'_I$. Expanding
\[
	\phi_k(e'_I)=\sum\limits_{J\in\II_k}c_{I,J}f'_J\, , \, c_{I,J}\in C
\]
we see, using that $C$ is saturated, that $(*)$ for our fixed $k$
is equivalent to 
\begin{equation}\tag{**} 
	\forall I,J\in\II_k: w(c_{I,J})\ge |J|-|I|,
\end{equation}
$w$ denoting the filtration of $C$. If $(G,\omega)$ is equi-$p$-valued, all the differences on the right-hand-side of $(**)$ are zero,
so that $(**)$ is trivially true.

By Proposition \ref{iso} any two such $\phi$ are chain homotopic via a
homotopy $h:\Sat X_\bullet\to \Sat Y_\bullet$ such that
$\gr(h)=0$. It remains
to check that it restricts to a homotopy 
$h:C\otimes_A X_\bullet\to C\otimes_B Y_\bullet$. We use the same generators 
as before. The additional assumption that $\gr(G)$ is generated in
degree $\frac{1}{e}$ implies $|I|=\frac{k}{e}$ for $I\in \II_k$.

Consider $e'_I$ for $I\in \II_k$. Then $h_k(e'_I)\in \Sat Y_{k+1}$ 
and expands as
\[
        h_k(e'_I)=\sum\limits_{J\in\II_{k+1}}d_{I,J}f'_J\, , \, d_{I,J}\in C
\]
$\gr(h)=0$, hence $\pi|d_{I,J}$ for all $I,J$.
As $e'_I=\pi^{-k}e_I$ and $f'_J=\pi^{-(k+1)}f_J$
this implies
\[
	h_k(e_I)=\sum\limits_{J\in \II_{k+1}}d_{I,J}\pi^{-1} f_J\ .
\] 
with $d_{I,J}\pi^{-1}\in C$ as required. 
\end{proof}

\begin{rem}\label{superfluous} It seems difficult to directly relate the complexes $C\otimes_A X_\bullet$ and $C\otimes_B Y_\bullet$ using the filtration techniques
successfully employed for example in \cite{ST} and \cite{T}, essentially
because these complexes do not satisfy any reasonable exactness properties.

In fact, we have $H_*(C\otimes_A X_\bullet)=\Tor^A_*(C,\Oh)$ and
$H_*(C\otimes_B Y_\bullet)=\Tor^B_*(C,\Oh)$ and one can check that,
unless $G=\{ e\}$, the algebra $C$ is not flat over neither $A$ nor $B$.

We have examples of saturated but not 
equi-$p$-valued groups and an isomorphism $\phi$ as above which 
does not restrict as in Proposition \ref{restrictiso}, but in 
all these examples it was possible by inspection to modify $\phi$ suitably.

It thus remains a tantalizing open problem to decide whether the
assumption ``equi-$p$-valued'' is superfluous in Proposition \ref{restrictiso}.
Of course, a positive answer would greatly extend the range of applicability of our integral Lazard comparison isomorphism.
\end{rem}

\begin{proof}[Proof of Theorem \ref{intlazard}] There is a filtration $\omega$ of $G$ such that $(G,\omega)$ is $p$-saturated, equi-$p$-valued of finite rank and $\omega(G)\subseteq\frac{1}{e}\Z$ for some integer $e\ge 1$, hence we are in the situation
studied in this subsection and in particular recall $\Oh,A,B,C,X_\bullet$ and
$Y_\bullet$ from above.  The continuous group cohomology $H^*_c(G,M)$ is defined using continuous cochains and the Bar-differential as in 
\cite[V 2.3.1.]{L}. By \cite[V 1.2.6 and 2.2.3.1]{L}  we have
\[ 
	H^*_c(G,M)\simeq H^*_c(\Z_p[[G]],M)\simeq \Ext^*_{\Z_p[[G]]}(\Z_p,M)
\]
and analogously
\[
	H^*_c(G,M\otimes_{\Z_p}\Oh)\simeq \Ext^*_A(\Oh, M\otimes_{\Z_p}\Oh)
\]
by the flatness of $\Oh$ over $\Z_p$. Introduce $N:=M\otimes_{Z_p}\Oh$.
Since $X_\bullet$ is a finite
free resolution of $\Oh$ over $A$, we obtain
\[
	\Ext_A^*(\Oh,N)=H^*\Hom_A(X_\bullet,N)=H^*\Hom_C(C\otimes_A X_\bullet,N)
\]
using that the $A$-module structure on $N$ extends to
a $C$-module structure, and then
\begin{multline*}
...\stackrel{Prop. \ref{restrictiso}}{\simeq} H^*\Hom_C(C\otimes_B Y_\bullet,N)=H^*\Hom_B(Y_\bullet,N)\\
\simeq H^*(\LL^*(G)\otimes_{\Z_p}\Oh,N)\simeq H^*(\LL^*(G),M)\otimes_{\Z_p}\Oh,
\end{multline*}
the last but one isomorphism by \cite[Lemma 9.2]{T}. Summing up we have

\begin{equation}\label{isoolin}
 H^*_c(G,M)\otimes_{\Z_p}\Oh\simeq H^*(\LL^*(G),M)\otimes_{\Z_p}\Oh
\end{equation}
of $\Oh$-modules.

We now turn to functoriality. Let $f:G\to H$ be filtered group homomorphism. 
We use $A(G)$,  $C(G)$, $X_\bullet(G)$, $Y_\bullet(G)$ (resp. $A(H)$, $C(H)$, $X_\bullet(G)$, $Y_\bullet(G)$) for the rings $A$, 
$C$ and complexes $X_\bullet$, $Y_\bullet$ corresponding to the group $G$ (resp. $H$).
The group homomorphism
induces a commutative diagram
\[\xymatrix{
 \gr(G)\ar[r]^{\gr(f)}\ar[d]^\isom&\gr(H)\ar[d]^\isom\\
 \gr(\LL^*(H))\ar[r]&\gr(\LL^*(H))
 }\]
As in Proposition \ref{iso} it lifts to a diagram
of filtered complexes of $\Sat A(G)$-modules
\[\xymatrix{
  \Sat X_\bullet(G)\ar[r]\ar[d]^\isom &\Sat X_\bullet(H)\ar[d]^\isom\\
  \Sat Y_\bullet(G)\ar[r]&\Sat Y_\bullet (H)
  }\]
which commutes up to homotopy and such that taking gradeds gives back the
previous diagram, and such that taking gradeds of the homotopy is $0$.
As in Proposition \ref{restrictiso} it restricts to a diagram of
filtered complexes of $C(G)$-modules
\[\xymatrix{
  C(G)\otimes X_\bullet(G)\ar[r]\ar[d]^\isom &C(G)\otimes X_\bullet(H)\ar[d]^\isom\\
  C(G)\otimes Y_\bullet(G)\ar[r]& C(G)\otimes Y_\bullet (H)
  }\]
 which commutes up to homotopy.

Compatibility with cup-products is the case $\Delta:G\to G\times G$.
Note that the generators of $G\times G$ are of the form
$(x,1)$ and $(1,y)$ for generators $x,y$ of $G$. Their filtration
is the same as that of $x,y$. 
\end{proof}
\section{The Lazard isomorphism for algebraic group schemes}\label{groupsection}
In this section we give, in the case of $p$-adic Lie-groups arising from
algebraic groups, a direct description of a map from analytic group 
cohomology to Lie-algebra cohomology by differentiating cochains.
In Proposition \ref{lazard} we recall that this coincides with Lazard's isomorphism.
This analytic description directly generalizes to $K$-Lie groups where $K$ is
a finite extension of $\Q_p$. In Theorem \ref{isom} we show that is also
an isomorphism in the case of $K$-Lie groups defined by algebraic groups.

\subsection{Group schemes}
Let $p$ be a prime number,
$K$ be a finite extension of $\Q_p$, let $R$ be its ring of integers with
prime element $\pi$. Throughout $\G$ will be a 
separated smooth group scheme over $R$
and $\g$ its $R$-Lie-algebra in the following sense:
\[ \g=\Lie(\G)=\Der_R(\Oh_{\G,e},R)\ .\]
Then $\g_K:=\g\otimes_RK=\Lie(\G_K)$ is its Lie-algebra as a $K$-manifold.

Note that this category is stable under base change and Weil restriction
for finite flat ring extensions $R\to S$. If $A\to B$ is a ring extension with
$B$ finite and locally free over $A$ and $X$ an $A$-scheme, we write
$X_B=X\times_A\Spec B$. If $Y$ is a $B$-scheme, we write
$\Res_{B/A}X$ for the Weil restriction, i.e.,
$\Res_{B/A}X (T)=X(T_B)$ for all $A$-schemes $T$.
See \cite[\S 7.6]{BLR} for properties of the
Weil restriction. In particular, if $\G$ is a group scheme over a discrete
valuation ring $R$, then $\G$ is quasi-projective by \cite[\S 6.4, Theorem 1]{BLR}.
This suffices to guarantee that $\Res_{S/R}(\G)$ exists for finite extensions $S/R$.

The following bit of algebraic geometry will be needed in the proofs.

\begin{lemma}\label{factor}Let $L/K$ be finite extension, $S$ the ring of integers of $L$. Consider
a separated smooth group scheme $\G$ over $S$. Then 
$\G$ is a direct factor of $\Res_{S/R}(\G)_S$.
\end{lemma}
\begin{proof}
Let $X$ be an $S$-scheme. For all $S$-schemes $T$ we describe 
$T$-valued points
of $\Res_{S/R}(X)_S$:
\begin{multline*}
 \Mor_S(T,\Res_{S/R}(X)_S)=\Mor_R(T,\Res_{S/R}(X))
=\Mor_S(T\times_R\Spec S,X)\\
=\Mor_S(T\times_S\Spec (S\tensor_R S),X)
\end{multline*} 
The natural map $\iota:S\to S\tensor_R S$ which maps $s\mapsto s\tensor 1$
induces the transformation of functors
\[ \Mor_S(T\times_S\Spec S,X)\xrightarrow{\iota} \Mor_S(T\times_S\Spec (S\tensor_R S),X)=
\Mor_S(T,\Res_{S/R}(X)_S)\]
and hence a morphism
\[ \iota:X\to \Res_{S/R}(X)_S\]
This is nothing but the adjunction morphism.

The multiplication $\mu_S:S\otimes_R S\to S$ is a section of $\iota$.
This again induces
a transformation of functors
\[ \Mor_S(T,\Res_{S/R}(X)_S)=\Mor_S(T\times_S\Spec(S\tensor_R S),X)\xrightarrow{\mu_S} \Mor_S(T\times_SS,X)\]
and hence a morphism
\[ \mu_S:\Res_{S/R}(X)_S\to X\ .\]
(Put $T=\Res_{S/R}(X)_S$ and the identity on the left.)
By construction $\mu_S$ is a section of $\iota$. Both are natural in $X$, 
hence $\G$ is a direct factor
of $\Res_{S/R}(\G)_S$ as group schemes.
\end{proof}

\begin{rem}If $L/K$ is Galois of degree $d$, then $\Res_{L/K}(\G)_L\isom \G^d$.
This carries over to the integral case if the extension is unramified.
The assertion becomes false for ramified covers. Note, however, that
the weaker statement of the lemma remains true.
\end{rem}
\subsection{Analytic description of the Lazard morphism}
Let $\G$ be a smooth connected  group scheme over $R$ with Lie-algebra $\g$. Let
$\cG\subset \G(R)$ be an open sub-Lie-group.

We denote
by $\Oh_\la(\cG)$ (locally) analytic functions on $\cG$, i.e. those
that can be locally written as a converging power series with coefficients in $K$.  
We denote by $H^i_\la(\cG,K)$ (locally) analytic group cohomology, i.e.,
cohomology of the bar complex $\Oh_\la(\cG^n)_{n\geq 0}$ with the usual differential. We denote $H^i(\g,K)$ Lie-algebra cohomology, i.e. cohomology of the
complex $\Lambda^*(\g_K^\vee)$ with differential induced by the dual of
the Lie-bracket.

\begin{defn}\label{analLaz}The {\em Lazard morphism} is the map 
\[ \Phi:H^i_\la(\cG,K)\to H^i(\g,K)\]
induced by the morphism of complexes
\begin{align*} \Oh_\la(\cG^n)&\to(\g_K^{n})^\vee\to \Lambda^n\g_K^\vee\\
                 f&\mapsto df_e
\end{align*}
\end{defn}

\begin{rem}It is not completely obvious that $\Phi$ is a morphism
of complexes. See \cite[Section 4.6. and Section 4.7.]{HK}.
\end{rem}
\begin{rem}$\Phi$ is compatible with the multiplicative structure.
\end{rem}

Recall from Lemma \ref{saturationlemma} 
that in the case $K=\Q_p$, the kernel $G$ of $\G(\Z_p)\to \G(\F_p)$
is filtered and has a subgroup $\cG$ of finite index which is
saturated and equi-$p$-valued. Indeed for $p\neq 2$, we have
$G=\cG$.

Let $\LL^*=\LL^*(\cG)$ be its integral Lazard Lie-algebra (see Definition \ref{defintlaz}). As reviewed in Example
\ref{standardlie} there is a natural isomorphism
\[ \g\otimes\Q_p\isom \LL^*\otimes\Q_p\ .\]
\begin{prop}\cite[Theorem 4.7.1]{HK} \label{lazard} For $K=\Q_p$ and $\cG$ saturated, the Lazard morphism $\Phi$ (see Definition \ref{analLaz}) 
agrees under the identification of Lie-algebras in Example \ref{standardlie} with the isomorphism defined by Lazard (\cite[V 2.4.9, V 2.4.10]{L})
\end{prop}
In particular, $\Phi$ is an isomorphism in this case.
\begin{remark}This is a case where our integral version of the result (Theorem \ref{intlazard}) can
be applied. As shown there this is again the same isomorphism.
\end{remark}

\subsection{The isomorphism over a general base}
\begin{thm}\label{isom}
Let $\G$ be a smooth group scheme over $R$ with connected generic fiber
and $\cG\subset \G(R)$ an open  subgroup. Then
the Lazard morphism 
$\Phi$ (see Definition \ref{analLaz}) is an isomorphism.
\end{thm}
The proof will take the rest of this note. 

\begin{rem}
Let us sketch the argument. We are first going
to show injectivity. For this we can restrict to smaller and smaller
subgroups $\cG$ and even to their limit. In the limit, the statement 
follows by base change from Lazard's result for $R=\Z_p$. We then show surjectivity.
Finite dimensionality of Lie-algebra cohomology implies that the
morphism is surjective for small enough $\cG$.  Algebraicity then
implies surjectivity also for the maximal $\cG$.
\end{rem}

By construction, the Lazard morphism $\Phi$ depends only on an infinitesimal
neighborhood of $e$ in $\cG$. Hence it factors through the Lazard
morphism for all open sub-Lie-groups of $\cG$ and even through its limit
\[ \Phi_\infty:\varinjlim_{\cG'\subset \cG} H^i_\la(\cG',K)\to H^i(\g,K)\]
\begin{lemma}The limit morphism
$\Phi_\infty$ is an isomorphism.
\end{lemma}
\begin{proof}For $K=\Q_p$ this holds by Proposition \ref{lazard} and
the work of Lazard \cite[V 2.4.9 and V 2.4.10]{L}.

The system of open sub-Lie-groups of $\cG$ is filtered, hence
\[ \varinjlim_{\cG'\subset \cG} H^i_\la(\cG',K)=H^i( \Oh_\la(\cG^\bullet)_e )\]
where $\Oh_\la(\cG^n)_e$ is the ring of germs of locally analytic functions
in $e$. Note that $\G(\Z_p)$ also carries the structure of a rigid
analytic variety and germs of locally analytic functions are nothing but
germs or rigid analytic functions. Hence they can be identified with a limit
of Tate algebras.

First suppose that $\G=\H_R$ for a smooth group scheme $\H$ over
$\Z_p$. 
Then
\[ \Oh_\la(\cG^n)_e\isom \Oh_\la(\H(\Z_p)^n)\tensor K\]
because Tate algebras are well-behaved under base change (see 
\cite[Chapter 6.1, Corollary 8]{BGR}).
Moreover, $\Phi_\infty$ is compatible with base change. As it
is an isomorphism for $\H$ it is also an isomorphism for $\G$.

Now consider general $\G$. By Lemma \ref{factor}  $\G$ is direct factor of
some group of the form $\H_R$ with $\H$ a group over $\Z_p$. Indeed,
 $\H=\Res_{R/\Z_p}(\G)$. By naturality, $\Phi_{\infty,\G}$  is a direct factor
 of $\Phi_{\infty,\H_R}$ and hence by the special case an isomorphism.
\end{proof}

\begin{cor}\label{inj}$\Phi$ is injective.
\end{cor}
\begin{proof}
$\G(R)$ is compact, hence all open sub-Lie-groups are of finite
index. If $\cG'\subset \cG$ is an open normal subgroup, we have
\[ H^i_\la(\cG,K)\isom H^i_\la(\cG',K)^{\cG/\cG'}\]
Hence the restriction maps
\[ H^i_\la(\cG,K)\to H^i_\la(\cG',K)\]
are injective. As the system of open normal subgroups is filtered, this
also implies that
\[ H^i_\la(\cG,K)\to \varinjlim_{\cG'}H^i_\la(\cG',K)\]
is injective.
The injectivity of $\Phi$ follows from the injectivity of $\Phi_\infty$.
\end{proof}

\begin{lemma}
Let $\cG\subset \G(R)$ be an open subgroup. Then there is an open subgroup $\cH\subset\cG$ 
 such that the Lazard morphism
for $\cH$ is bijective.
\end{lemma}
\begin{rem} Note that this is precisely what Lazard proves  over $\Q_p$
with $\cH$ the saturated subgroup of $\cG=\G(\Z_p)$.
\end{rem}
\begin{proof}
As $\Phi_\infty$ is bijective and $\Phi$ injective, it suffices to show that 
there is
$\cH$ such that the restriction map
\[ H^i_\la(\cH,K)\to   \varinjlim_{\cG'}H^i_\la(\cG',K)\]
is surjective. Let $\alpha$ be a cocycle with class 
$[\alpha]\in \varinjlim_{\cG'}H^i_\la(\cG',K)$. 
By definition it is represented by a cochain on some
$\cG'$. It is a cocycle (possibly on some smaller $\cG'$). Hence
$[\alpha]$ is in the image of the restriction map for $\cG'$.

Lie-algebra cohomology is finite dimensional by definition, hence this
is also true for $\varinjlim_{\cG'}H^i_\la(\cG',K)$. By intersecting the
$\cG'$ for a basis we get the group $\cH$ we wanted to construct.
\end{proof}

\begin{proof}[Proof of Theorem \ref{isom}.]
Injectivity has already been proved in Corollary \ref{inj}.
 We use
an argument of Casselman-Wigner \cite[\S 3]{CW}  to conclude.
The operation of
$\G_K$ on $H^i(\g,K)$ is algebraic. Hence the stabilizer $\bbS_K$ is
a closed subgroup of $\G_K$. On the other hand it contains some open
subgroup of $\G(R)$. This implies that $\bbS_K=\G_K$ because $\G_K$ is connected. 
Hence $\cG\subset \G(R)$ operates trivially and thus
\[ \Phi:H^i_\la(\cG,K) \to H^i(\g,K)\]
is surjective.
\end{proof}
\begin{rem}The argument also works for cohomology with coefficients
in a finite dimensional algebraic representation of the group.
\end{rem}

\setcounter{section}{0}


\begin{appendix}
\section{}
In this short Appendix we explain an interesting remark of T. Weigel communicated
to us during the proof reading of the present paper.

We fix an {\em odd} prime $p$ and only remark without going into details 
that at the prime $2$ similar but strictly weaker results are available.

\begin{appprop}[Weigel]
Let $(G,\omega)$ be a $p$-saturated group of finite rank which is equi-$p$-valued. 
\begin{itemize}
\item[i)] The filtration of $G$ defined by $\omega$ is the lower $p$-series.
\item[ii)] There is a valuation $\omega'$ on $G$ with all properties of $\omega$ and
such that in addition there is an ordered basis $g_i\in G$ satisfying $\omega'(g_i)=1$
for all $i$.
\item[iii)] The group $G$ is uniformly powerful.
\end{itemize}
\end{appprop}
\begin{proof}
i) Denote by $t$ the valuation of an ordered basis of $G$.
Since $(G,\omega)$ is equi-$p$-valued we have
\[
	\gr(G)= \F_p[\epsilon]\cdot \gr^t(G),
\]
i.e. the filtration of $G$ defined by
$\omega$ has jumps exactly as follows
\[
	G=G^t\supseteq G^{t+1}\supseteq G^{t+2}\ldots.
\]
For every $n\ge 0$ we see that 
\[ 
	G^{t+n}/G^{t+n+1}\simeq \gr^{t+n}(G)\simeq \epsilon^n\cdot\gr^t(G)
\]
is elementary $p$-abelian, hence the Frattini subgroup $\Phi(G^{t+n})$ of 
$G^{t+n}$ satisfies
\[ 
	\Phi(G^{t+n})\subseteq G^{t+n+1}.
\]
On the other hand we have inclusions 
\[ 
	G^{t+n+1}\subseteq(G^{t+n})^p\subseteq\Phi(G^{t+n}),
\]
the first one by divisibility and the second one as a general property of
Frattini subgroups, thereby proving i).\\

ii) Since $(G,\omega)$ is $p$-saturated we have
\[ 
\frac{1}{p-1}<t\leq\frac{p}{p-1}.
\]
We now check that $\omega':=\omega+1-t$ has the desired properties:\\
Write $c:=1-t$ and first assume $t\ge 1$, i.e. $c\leq 0$. Then for all $x,y\in G$
\[ 
\omega'([x,y])=\omega([x,y])+c\ge\omega(x)+\omega(y)+c=\omega'(x)+\omega'(y)-c\geq\omega'(x)+\omega'(y).
\]
This shows that $\omega'$ has properties 1)-3) of Definition \ref{filtrationdefn}. 
Since (recall $p\neq 2$)
\[ 
	\omega'\geq 1> \frac{1}{p-1},
\]
$\omega'$ has property 4) of loc. cit. and property 5) is trivially satisfied. 
Assume $x\in G$ with 
\[ 
	\omega(x)+c=\omega'(x)>\frac{p}{p-1}.
\]
Then again because $c\leq 0$
\[ 
	\omega(x)>\frac{p}{p-1}-c\geq\frac{p}{p-1},
\]
and $x\in G^p$, since $\omega$ is saturated. Hence $\omega'$ has property 6) of loc. cit.
and property 7) is trivially satisfied.

Now assume that $t<1$, i.e. $c>0$. For all $x$,$y$ in an ordered basis of $G$ we have
\[ \omega([x,y])\geq \omega(x)+\omega(y)=2t>t\]
which, using the structure of $\gr(G)$, implies the first inequality in
\[ \omega([x,y])\geq t+1=\omega(x)+\omega(y)+1-t=\omega(x)+\omega(y)+c.\]
By uniformity, this inequality holds for all $x,y\in G$ and thus $\omega'$ has property 2), 
hence properties 1)-3) of Definition \ref{filtrationdefn}. Assume $x\in G$ with 
\[ 
	\omega(x)+c=\omega'(x)>\frac{p}{p-1}.
\]
Then 
\[ 
	\omega(x)>\frac{p}{p-1}-c=\frac{p}{p-1}+t-1=t+\frac{1}{p-1}>t.
\]
By the structure of $\gr(G)$ this implies
\[
	\omega(x)\geq t+1.
\]
Using that $(G,\omega)$ is $p$-saturated we have $t+1\geq \frac{p}{p-1}$ and $x\in G^p$.
This implies that $\omega'$
has property 6) of loc. cit. Since $\omega'\geq 1$, $\omega'$ has property 4) and
the rest of ii) clear.\\

iii) We first show that $G$ is uniform.
Choose $\omega'$ as in part ii), then 
\[ 
	G=G^1\supseteq G^2\supseteq\ldots
\]
is the lower $p$-series and 
\[ 
	[G^n:G^{n+1}]=|\epsilon^n\cdot \gr^1(G)|=[G^1:G^2].
\]
Furthermore, $G/G^p$ is abelian because
\[ 
	[G,G]\subseteq G^2\subseteq (G^1)^p= G^p.
\]
\end{proof}
We conclude by explaining the simplifications implied by this
result for our integral Lazard isomorphism, i.e. Theorem \ref{intlazard}:

\begin{appthm}Let $p$ be an odd prime, $(G,\omega)$ a compact, saturated, equi-$p$-valued group. Let $M$ be as in Theorem \ref{intlazard}. Then
there is an isomorphism of graded $\Z_p$-modules
\[ \phi_G(M):H^*_c(G,M)\simeq H^*(\LL^*(G),M).\]
It is natural in $M$ and $G$ and compatible with cup-products..
\end{appthm}
\begin{proof}
Replace $\omega$ by the filtration constructed above. Then we
are in the case $e=1$ of Theorem \ref{intlazard} with generators in degree $1$. 
For naturality
in $G$ note that
every group homomorphism $f$ (is continuous and) respects the lower $p$-series, hence the isomorphism is natural with respect to $f$.
In particular then, the isomorphism is always compatible with 
cup-products.
\end{proof}

\end{appendix}

\end{document}